\documentclass{amsart}
\usepackage{amssymb}
\newtheorem{thm}{Theorem}[section]
\newtheorem{cor}[thm]{Corollary}

\newtheorem{lem}[thm]{Lemma}

\theoremstyle{definition}

\theoremstyle{remark}

\numberwithin{equation}{section}

\begin{document}
\title[1-factorizations of Cayley graphs]{1-factorizations of Cayley graphs}%
\author{A. Abdollahi}%
\address{Department of Mathematics, University of Isfahan, Isfahan 81746-71441, Iran;
and Institute for Studies in Theoretical Physics and Mathematics
(IPM); Tehran, Iran.}%
\email{a.abdollahi@math.ui.ac.ir}%
\thanks{This research  was in part supported by a grant from
IPM (No. 85200032). The author thanks the Center of Excellence
for Mathematics, University of
Isfahan.}%

\begin{abstract}
In this note we prove that all connected Cayley graphs of every
finite group $Q \times H$ are 1-factorizable, where $Q$ is any
non-trivial group of 2-power order and $H$ is any group of odd
order.
\end{abstract}
\maketitle
\section{\bf Introduction and Results}
Let $G$ be a non-trivial group, $S\subseteq G\backslash\{1\}$ and
$S^{-1}=\{s^{-1} : s\in S\}$. The Cayley graph $\Gamma(S:G)$ of
the group $G$ with respect to the set $S$ has the vertex set $G$
and the edge set $\big\{\{g,sg\} : g\in G,
s\in S \cup S^{-1}\big\}$.\\

A $j$-factor of a graph is a spanning subgraph which is regular
of valence $j$. In particular, a 1-factor of a graph is a
collection of edges such that each vertex is incident with
exactly one edge.  A 1-factorization of a regular graph is a
partition of the edge set of the graph into disjoint 1-factors. A
1-factorization of a regular graph of valence $v$ is equivalent
to a coloring of the edges in $v$ colors (coloring each 1-factor
a different color). This enables us to use a very helpful result:
Any simple, regular graph of valence $v$ can be edge-colored in
either $v$ or $v+1$ colors.
This is a specific case of Vizing's theorem (see \cite[pp. 245-248]{Ore}). \\

We study the conjecture that says all Cayley graphs $\Gamma(S:G)$
of groups $G$ of even order are 1-factorizable whenever
$G=\left<S\right>$. There are some partial results on this
conjecture obtained by Stong \cite{Stong}. Here we prove\\

\noindent{\bf Theorem.} {\sl Let $H$ be a finite group of odd
order and let  $Q$ be a finite group of order $2^k$ ($k>0$). Then
the Cayley graph $\Gamma(S:Q \times H)$ is 1-factorizable for all
generating sets $S$ of $Q \times H$.}\\

As a corollary we prove that  all connected Cayley graphs of every
finite nilpotent group of even order are 1-factorizable which has
been proved by Stong in \cite[Corollary 2.4.1]{Stong} only for
Cayley graphs on minimal generating sets.

\section{\bf Proof of the Theorem}
We need the following lemma whose proof is more or less as Lemma
2.1 of \cite{Stong} with some modifications.
\begin{lem}\label{l1}
Let $H$ be a finite group of odd order. Then the Cayley graph
$\Gamma(S:\mathbb{Z}_2\times H)$ is $1$-factorizable, for any
generating set $S$ of $\mathbb{Z}_2 \times H$ containing  exactly
one element of even order.
\end{lem}
\begin{proof}
Let $a$ be the only element of $S$ of even order. Then $a=zh$,
where $z\in \mathbb{Z}_2$ and $h\in H$ and $z$ of order $2$. If
$a^2=1$, then $h=1$ and $S\backslash\{a\}\subseteq H$ and so
$axa^{-1}=x$ for all $x\in S\cap H$. Thus, in this case, Theorem
2.3 of \cite{Stong} completes the proof. Therefore we may assume
that $a^2\not=1$.  Let
$\Gamma'=\Gamma(S\backslash\{a\}:\mathbb{Z}_2\times H)$ and
$\Gamma_1$ and $\Gamma_2$ be the induced subgraphs of $\Gamma'$
on the sets $H$ and $zH$, respectively. It can be easily seen
that the map $x\mapsto zx$ is an graph isomorphism from
$\Gamma_1$ to $\Gamma_2$. By Vizing's theorem the edges in both
$\Gamma_1$ and $\Gamma_2$ can be edges-colored in the same manner
in $|S\backslash\{a\}|+1$ colors (by ``the same manner'' we mean
that the edge $\{h_1,h_2\}$ in $\Gamma_1$ has ``the same'' color
as $\{zh_1,zh_2\}$ in $\Gamma_2$, and vice versa). Then all that
remains to be done is to color the edges from $H$ to $zH$, that
is the following two `disjoint' $1$-factors of
$\Gamma(S:\mathbb{Z}_2\times H)$ (here we use $a^2\not=1$):
$$\big\{\{x,ax\}\;|\; x\in H\big\} \;\;\text{and}\;\;
\big\{\{x,a^{-1}x\}\;|\; x\in H\big\}. \eqno{(*)}$$ (note that the
edges of $\Gamma(S:\mathbb{Z}_2\times H)$ are exactly the edges
of $\Gamma_1$, $\Gamma_2$ and those in the above $1$-factors).
Now since both $x\in H$ and $zx\in zH$ have edges (in $\Gamma_1$
and $\Gamma_2$, respectively) of the same $|S\backslash\{a\}|$
colors to them, there are `two' colors (note that here we again
use $a^2\not=1$) that can be used to color $1$-factors in $(*)$.
This completes the proof.
\end{proof}

 \noindent{\bf Proof of the Theorem.} Let $G=Q\times H$ and  $S$ be
any generating set of $G$. We argue by induction on $|S|$. If
$|S|=1$, then $G$ is a cyclic group of even order and Corollary
2.3.1 of \cite{Stong} completes the proof. Now assume that
$|S|>1$ and for any  non-trivial group $Q_1$ of 2-power order and
subgroup $H_1$ of $H$ the Cayley graph  $\Gamma(S_1:Q_1 \times
H_1)$ is 1-factorizable for any generating set $S_1$ of
$Q_1\times H_1$ with  $|S_1|<|S|$.
 Since the set of elements of odd order in $G$ is the subgroup $H$ and $G=\left<S\right>$, $S$ has at
least one element $a$ of even order. First assume  that $S$ has
another element distinct from $a$ of even order. Consider  the
subgroup $G_1$ generated by $S\backslash\{a\}$  of $G$. Then
$G_1=Q_1 \times H_1$ for some subgroups $Q_1\leq Q$ and $H_1\leq
H$ such that $Q_1\not=1$.  Therefore the induction hypothesis
implies that $\Gamma(S\backslash\{a\}:G_1)$ has a
1-factorization. Since $\Gamma(S\backslash\{a\},G)$ consists of
disjoint copies of $\Gamma(S\backslash\{a\}:G_1)$ which are
1-factorizable, $\Gamma(S\backslash\{a\},G)$ has a
1-factorization. Now since the only element of $S\backslash
(S\backslash\{a\})$ has even order, Lemma 2.2 of
\cite{Stong} shows that $\Gamma(S:G)$ is 1-factorizable.\\
 Hence we may assume that $a$ is the only element of $S$ of even order. Since $a=a_1a_2$ for some $a_1\in Q$ and $a_2\in H$,
 we have
$$G=\left<S\right>=\left<S\backslash\{a\},a_1a_2\right>=\left<a_1\right>
\times \left<S\backslash\{a\},a_2\right>.$$ It follows that
$Q=\left<a_1\right>$. Consider the subgroup
$N=\left<a_1^2\right>$. Then $N$ is a normal subgroup of $G$ such
that $N\cap S=\varnothing$. It is easy to see that when $s,t\in
S$ with $s\not=t^{\pm 1}$, neither $st$ nor $st^{-1}$ belongs to
$N$.  Now by Lemma 2.4 of \cite{Stong}, it is enough to show that
$\Gamma(\frac{SN}{N}:\frac{G}{N})$ is 1-factorizable. Since
$\frac{G}{N}\cong \mathbb{Z}_2 \times H$, it follows from Lemma
\ref{l1} that $\Gamma(\frac{SN}{N}:\frac{G}{N})$ is
1-factorizable.  This completes the proof. $\hfill\square$\\

\begin{cor}
If $G$ is a finite nilpotent group of even order, then
$\Gamma(S:G)$ is $1$-factorizable for all generating sets $S$ of
$G$.
\end{cor}
\begin{proof}
It follows from the Theorem and the fact that every finite
nilpotent group is the direct product of its Sylow subgroups.
\end{proof}


\end{document}